\documentclass[11pt,a4paper]{article}
\usepackage{amsfonts}
\usepackage{amsmath}
\usepackage{amssymb}
\usepackage[english]{babel}
\usepackage[T1]{fontenc}
\usepackage{graphicx}
\usepackage{subfigure}
\usepackage{epsfig}
\usepackage{latexsym}
\usepackage{amstext}
\usepackage{amsthm}
\usepackage{graphics}
\usepackage{enumerate}
\usepackage{paralist}
\usepackage{color}
\usepackage{fullpage}
\usepackage{dsfont}
\newtheorem{prop}{Proposition}[section]
\newtheorem{rem}[prop]{Remark}
\newtheorem{lem}[prop]{Lemma}
\newtheorem{theo}[prop]{Theorem}
\newtheorem{cor}[prop]{Corollary}

\numberwithin{equation}{section}

\newcommand{\beq}{\begin{eqnarray}}
\newcommand{\beqq}{\begin{eqnarray*}}
\newcommand{\eeq}{\end{eqnarray}}
\newcommand{\eeqq}{\end{eqnarray*}}

\def\QED{\quad\hbox{\hskip 4pt\vrule width 5pt height 6pt depth 1.5pt}}

\title{Windings of planar processes, \\ Exponential Functionals and Asian options}
\author{{\sc Wissem Jedidi} \thanks{Department of Statistics \& OR, King Saud University, P.O. Box 2455, Riyadh 11451, Saudi Arabia and
Universit\'e de Tunis El Manar, Facult\'e des Sciences de Tunis, LR11ES11 Laboratoire d'Analyse Math\'ematiques et Applications,
2092, Tunis, Tunisia. E-mail: wissem\_jedidi@yahoo.fr}
\ {\sc and} \ {\sc Stavros Vakeroudis} \thanks{Corresponding Author} \thanks{Department of Mathematics,
Track: Statistics and Actuarial-Financial Mathematics,
University of the Aegean, Vourliotis Building, Office: Y5,
   83200 Karlovasi, Samos,
    Greece. E-mail: stavros.vakeroudis@gmail.com \ \ Web: https://svakeroudis.wordpress.com  }  }
\date{\today}
\begin{document}

\maketitle
\begin{abstract}
Motivated by a common Mathematical Finance topic, we discuss the reciprocal of the exit time from a cone of planar Brownian motion
which also corresponds to the exponential functional of an associated Brownian motion. 
We prove a conjecture in \cite{VaY12b} concerning infinite divisibility properties of this random variable and we present
a novel simple proof of De Blassie's result in \cite{DeB87,DeB88} about the asymptotic behaviour of the distribution of the Bessel clock
appearing in the skew-product representation of planar Brownian motion, as $t\rightarrow\infty$.
Similar issues for the exponential functional of a L\'evy process are also discussed.
We finally use the findings obtained by the windings approach in order to get results for quantities associated to the pricing of Asian options.
\end{abstract}

$\vspace{5pt}$

$\vspace{5pt}$

\section{Introduction}
Windings of 2-dimensional processes, and especially of planar Brownian motion have several applications in Financial Mathematics for instance, where the exponential functionals of Brownian motion are of special interest.
A fundamental example is the pricing of Asian options (see e.g. \cite{Yor93,Duf00,GeY01,Yor01,MaY05,MaY05a}), where the payout of an Asian call option is given by:
\beqq
E\left[\left(\frac{1}{t}\int^{t}_{0} ds \ \exp(\beta_{s}+\nu s)-K\right)^{+}\right] \, ,
\eeqq
where $\left(\beta_{u},u\geq0\right)$ is a real Brownian motion, $\nu\in\mathbb{R}$ and the non-negative number $K$ is the strike price.
It is easy to show (for further details, see e.g. \cite{GeY01}) that the computation of this expectation simplifies to the computation of
\beqq
E\left[\left(\int^{t}_{0} ds \ \exp(\beta_{s}+\nu s)-K\right)^{+}\right] \, ,
\eeqq
which follows by studying the quantity
\beqq
E\left[\int^{t}_{0} ds \ \exp(\beta_{s}+\nu s)\right] \,.
\eeqq
In particular, in \cite{Yor93} one can find a more detailed discussion for the distribution of the exponential functional
$$A_t^{(\nu)}:=\int^{t}_{0} ds \ \exp(\beta_{s}+\nu s)$$ taken up to
a random time $T_\lambda$ which follows the exponential distribution with parameter $\lambda>0$ and is independent from $\beta$.
More precisely, Yor \cite{Yor92} obtained that
$$2A_{T_\lambda}^{(\nu)}\stackrel{(law)}{=}\frac{Q_{1,a}}{2G_{b}}\stackrel{(law)}{=}\frac{1-\mathcal{U}^{1/a}}{2G_{b}},$$
where $Q_{1,a}\sim\mathrm{Beta}(1,a), \ G_{b}\sim\mathrm{Gamma}(b), \ \mathcal{U}\sim U[0,1]$,
$a=(\nu/2)+(1/2)\sqrt{2\lambda+\nu^2}$, $b=a-\nu$ and the random variables in the identities in law are assumed to be independent. The class of  Generalized Gamma Convolution distributions  ($GGC$) is an important subclass of infinitely divisible distributions. The monograph of Bondesson \cite{Bon92} is devoted to the deep study of this class.  At this stage, we could observe that $G_b \sim GGC$ and that $1/Q_{1,a}\sim GGC$ (see \cite{simon}) for the last result. We easily conclude by the stability property by product (see \cite{Bon15}), that $1/A_{T_\lambda}^{(\nu)} \sim GGC$. This result gives a flavor of what we will obtain in Subsection \ref{s2}.
\\
We  mainly deal with such exponential functionals and we try to investigate their distribution properties.
The key observation of our approach is the fact that exponential functionals are strongly related to the windings of associated processes (see for instance Proposition \ref{relations} below) and this offers a possible direction to describe them. Indeed, instead of using an independent random time such as $T_\lambda$ above, one may study the exponential functional up to the first hitting time of a specific level by another independent real Brownian motion, and this is related to the windings of planar Brownian motion as we shall see later on.
Once obtaining results of the exponential functional in the framework of windings, they may be used to study the asymptotic behaviour of the exponential functional of interest.
The goal of this paper is first to explore the distribution of exponential functionals in terms of planar processes and then to take profit of this in order to discuss the exponential that we meet in the pricing of Asian options.
Note that windings of different types of processes (e.g. jump processes) are related to different kinds of exponential functionals (e.g. exponential functionals of L\'{e}vy processes) hence, to different types of Asian options.

We shall first study exponential functionals in terms of planar Brownian motion, i.e. taken up to an independent random time different from the case mentioned above, that is (we suppose here that $\nu=0$ but at the end of the paper we will also discuss the case where $\nu\neq0$):
\beqq
\int^{T^{\gamma}_{c}}_{0} \exp(2\beta_{s}) ds \ ,
\eeqq
where $\left(\gamma_{u},u\geq0\right)$ is another real Brownian motion independent from $\beta$ and the exit time $T^{\gamma}$  is given by (\ref{exit}) below. For more precise connection with the windings, see Proposition \ref{relations}.

We shall also discuss the exponential functional associated to jump processes, i.e.
$$\int^{t}_{0}\exp(\alpha\xi_{s})ds,$$
where $\left(\xi_{u},u\geq0\right)$ is a (non-symmetric) L\'{e}vy process and $\alpha\in(0,2)$ is a constant.
As we shall see, this exponential functional is again related to the study of an associated planar Stable process with index of stability $\alpha$.
\\
We consider the following processes:
\begin{itemize}
\item $Z=(Z_{t},t\geq0)$ a planar Brownian motion (BM),
\item $U=(U_{t},t\geq0)$ a planar Stable process of index $\alpha\in(0,2)$,
\end{itemize}
both starting from a point different from 0 and, without loss of generality we may consider that they are both issued from 1.
For $W=Z$ or $U$ and $\alpha \in (0,2]$, we formally define the clock, its inverse and  the winding process, given respectively, for all $t\geq0$, by
\begin{equation}\label{clock}
H^{W}_{t}=\int^{t}_{0}\frac{ds}{\left|W_{s}\right|^{\alpha}},\quad A^{W}(u)=\inf\{t\geq0,\; H^{W}(t)>u\}  \quad \mbox{and}\quad \theta^{W}_{t}=\mathrm{Im}\big(\int^{t}_{0}\frac{dW_{s}}{W_{s}}\big)\, .
\end{equation}
Note that the planar Brownian motion case, that is $W=Z$, corresponds to $\alpha=2$.
We also define the exit times from a cone of single and of double border of a process $V$ ($V$ will be $W$ or a functional of $W$ in the sequel) by
\begin{equation}\label{exit}
T^{V}_{c} = \inf\{ t\geq 0:\;V_{t}\geq c \} \quad \mbox{and} \quad T^{|V|}_{c} =\inf\{t\geq 0:\;|V_{t}|\geq c \}, \quad c>0.
\end{equation}

For the planar Brownian motion case ($W=Z$), it is well-known  (see \cite{ItMK65} for instance), that since $Z_{0}\neq0$, the process $Z$ does not visit a.s. the point $0$ but keeps winding around it infinitely often. Hence, its continuous winding process $\theta^{Z}_{t}$ is well defined. We recall as well the skew-product representation of planar BM (see e.g. \cite{ReY99}):
\beq\label{skew-product}
\log\left|Z_{t}\right|+i\theta^{Z}_{t}=\int^{t}_{0}\frac{dZ_{s}}{Z_{s}}=\left(\beta_{u}+i\gamma_{u}\right) \Bigm|_{u=H^{Z}_{t}} \ ,
\eeq
where $(\beta_{u}+i\gamma_{u},u\geq0)$ denotes another planar Brownian motion starting from $\log 1+i0=0$ and $H^Z$ is given by \eqref{clock} with $\alpha=2$ (for the Bessel clock $H^{Z}$ see also \cite{Yor80}).  It is straightforward that the two $\sigma$-fields $\sigma \{\left|Z_{t}\right|,t\geq0\}$ and $\sigma \{\beta_{u},u\geq0\}$ are identical, whereas $(\gamma_{u},u\geq0)$ is independent from $(\left|Z_{t}\right|,t\geq0)$. Note that the inverse of $H^{Z}$ is represented by the functional
\beq\label{inverse}
A^{Z}(t)=\inf\{u\geq0, \;H^{Z}(u)>t\}=\int^{t}_{0}ds \ \exp(2\beta_{s}) \ .
\eeq
We address the interested reader e.g. to \cite{PiY86} for more details about planar Brownian motion.

Contrary to planar Brownian motion, one cannot define the winding number directly for the isotropic planar Stable process $U$ (see \cite{BeW96,DoV12,KyV16}). However, we can consider a path on a finite time
interval $[0,t]$ and "fill in" the gaps with line segments in order to obtain the curve of a continuous function
$f:[0,1]\rightarrow\mathbb{C}$ such that $f(0)=1$. The origin 0 is polar and  $U$  has no jumps across 0 a.s., thus we have $f(u)\neq0$
for every $u\in[0,1]$ and the process of the winding number of $U$ around 0, $\theta^{U}=(\theta^{U}_{t},t\geq0)$,
is well-defined, has c\`{a}dl\`{a}g paths of absolute length greater than $\pi$ and is given by
\beqq
\exp(i\theta^{U}_{t})=\frac{U_{t}}{|U_{t}|}, \quad  t\geq0  \ .
\eeqq

For $W=Z$ or $U$, we shall study the exit times from a cone of single and of double border, i.e. the stopping times $T^{\theta^{W}}$
and $T^{|\theta^{W}|}$ given by (\ref{exit}), and  also the asymptotic behaviour of the associated winding process.

The rest of the paper is organized as follows: we start by discussing windings and the associated version of Spitzer's Asymptotic Theorem (that corresponds to the large time asymptotics) (i) for planar Brownian motion, and (ii) for isotropic planar Stable processes
(note that for the latter it is not exactly an analogue of Spitzer's Asymptotic Theorem but mostly a large time asymptotics result).
In Section \ref{BM}, we characterize  the distribution of the exit times from a single and from a double border cone
which corresponds also to the exponential functional of Brownian motion in the framework of planar Brownian motion.
Then, we turn our interest to infinite divisibility properties of this quantity and we shall prove a conjecture in
\cite[Remark 3.2]{VaY12b}.
We shall also present in Subsection \ref{DeBlassie} a new simple proof of De Blassie's result in \cite{DeB87,DeB88}
stating that if $R(s)$ denotes a Bessel process, then, for every $u>0$ and for every $\lambda>0$,
$$P\left(\int^{t}_{0}R_{s}^{-2}ds\leq u\right)=O\left(t^{-\lambda}\right),\quad \mbox{as}\; t\rightarrow\infty.$$
Recall that
$$H^Z_t:=\int^{t}_{0}\frac{1}{|Z_{s}|^{2}}ds=\int^{t}_{0}R_{s}^{-2}ds,$$
hence the last  result corresponds to the asymptotic behaviour of the Bessel clock associated to planar Brownian motion $H^Z$, as $t\rightarrow\infty$.
The initial proof due to R. Dante De Blassie used results of Burkholder together with a theorem taken by the book of Port and Stone \cite{PoS78}.
Here, we propose a novel elementary self-contained proof.

Section \ref{St} focuses on the windings of isotropic planar Stable processes where a large time asymptotics result due to Bertoin and Werner \cite{BeW96} is presented for the sake of completeness in order to use it in the following section. Finally, Section \ref{Asian} deals with applications of the previous results to the pricing of Asian options.
More precisely, we discuss separately the case of exponential functionals of Brownian motion and the one of L\'{e}vy processes.

To recapitulate, the main results of the paper are the following:
\begin{itemize}
  \item we prove a conjecture in \cite{VaY12b} concerning infinite divisibility properties the inverse of the exponential functional in terms of planar Brownian motion;
  \item we propose a novel simple self-contained proof of De Blassie's result in \cite{DeB87,DeB88} concerning the distribution of the Bessel clock appearing in the skew-product representation of planar Brownian motion, for $t\rightarrow\infty$;
  \item we use results concerning exponential functionals in terms of windings in order to study exponential functionals needed for the pricing of Asian options by invoking William's "pinching method".
\end{itemize}

The approach proposed in this manuscript, has intrinsic theoretical interest since it provides further characterisation of the exponential functional associated to different Stochastic processes, including analytic properties (e.g. infinite divisibility properties).
On the other hand, it gives a new direction and perspectives in order to proceed to the pricing of different types of Asian options, such as the ones associated to jump (L\'{e}vy) processes.

\section{Planar Brownian motion}\label{BM}
\subsection{Windings and exponential functionals}
We first recall our main tool, which is Bougerol's celebrated identity in law \cite{Bou83}. It states that
if $(\beta_{u},u\geq0)$ and $(\hat{\beta}_{u},u\geq0)$ are two independent linear Brownian
motions both started from 0, then we have the identity
\beq\label{bougerol}
\sinh(\beta_{t}) \stackrel{(law)}{=} \hat{\beta}_{A^{Z}_{t}(\beta)=\int^{t}_{0}ds\exp(2\beta_{s})}, \quad \mbox{for every fixed $t\geq0$} \ .\eeq
For the proof and other developments of this identity, see \cite{Vak12} and the references therein. We will study Bougerol's identity in law in terms of planar Brownian motion, which is strongly related to exponential functionals of Brownian motion as one can see below.
To that end, we recall that the  exit times $T^{\gamma}_{c}$ and $ T^{|\gamma|}_{c}$  for the BM $\gamma$ associated to $\theta^{Z}$ are given by (\ref{exit}). As a first result, we obtain the following (see also \cite{Vak11th,Vak11}).
\begin{prop}\label{relations}
It holds that:
\beqq
T^{\theta^{Z}}_{c}=A^{Z}_{T^{\gamma}_{c}} \quad  \mathrm{and}  \quad  T^{|\theta^{Z}|}_{c}=A^{Z}_{T^{|\gamma|}_{c}}.
\eeqq
\end{prop}
{\noindent \textbf{Proof of Proposition \ref{relations}.}}
It follows by the skew-product representation $(\theta^{Z}_{t}=\gamma_{H^{Z}_{t}})$, using the fact that
$A^{Z}$ is the inverse of $H^{Z}$ (see also \eqref{inverse}), i.e.
\beqq
T^{\theta^{Z}}_{c}=\inf\{ t:\theta^{Z}_{t}=c \}=\inf\{ t:\gamma_{H^{Z}_{t}}=c \}
\stackrel{s=H^{Z}_{t}}{=} \inf\{ A^{Z}_{s}:\gamma_{s}=c \}=A^{Z}_{T^{\gamma}_{c}} \ .
\eeqq
The second relation follows similarly.
\hfill\QED
\\
From now on, all the results may be stated either for $A^{Z}_{T^{\gamma}_{c}}$ (resp. $A^{Z}_{T^{|\gamma|}_{c}}$) or for $T^{\theta^{Z}}_{c}$ (resp. $T^{|\theta^{Z}|}_{c}$).
For the sake of applications in the Mathematical Finance framework, we will mostly use the first notation.
\\
We recall Spitzer's celebrated asymptotic Theorem for planar BM \cite{Spi58}. For other proofs see e.g. \cite{VaY12} and the references therein.
\begin{theo}[Spitzer's Asymptotic Theorem (1958)]\label{Spitzertheo}
The following convergence in law holds:
\beqq
 \frac{2}{\log t} \; \theta^{Z}_{t} \overset{{(law)}}{\underset{t\rightarrow\infty}\longrightarrow} C_{1} \ ,
\eeqq
where $C_{1}$ denotes a standard Cauchy distributed random variable.
\end{theo}
{\noindent We introduce now the function}
\begin{equation}
\varphi(x) = \arg \sinh ^2 (\sqrt{x})=\log^2(\sqrt{x}+\sqrt{1+x}),\label{p3}
\end{equation}
which plays a key role in the rest of the paper. The next proposition comes essentially from \cite{Vak11th,Vak11} and we give its proof for the sake of completeness.
\begin{prop}\label{prop1}
The distributions of \ $A^{Z}_{T^{\gamma}_{c}}$ and of \ $A^{Z}_{T^{|\gamma|}_{c}}$ are characterized by the following Gauss-Laplace transforms: for all $x \geq 0$ and $m=\frac{\pi}{2c}$, we have
\begin{eqnarray}
c \; E \left[ \sqrt{\frac{\pi}{2 A^{Z}_{T^{\gamma}_{c}}}} \exp \left( -\frac{x}{2A^{Z}_{T^{\gamma}_{c}}} \right) \right] &=& \frac{1}{\sqrt{1+x}} \: \frac{c^{2}}{c^{2}+\varphi(x)} \ , \label{GLT1} \\
c \; E \left[ \sqrt{\frac{2}{\pi A^{Z}_{T^{|\gamma|}_{c}}}} \exp \left( -\frac{x}{2A^{Z}_{T^{|\gamma|}_{c}}} \right) \right] &=&  \frac{1}{\sqrt{1+x}} \, f_m(x),\label{GLT2}
\end{eqnarray}
where $\varphi$ is given by (\ref{p3}) and
\begin{equation}
f_m(x)=\frac{2}{(\sqrt{1+x}+\sqrt{x})^{m}+(\sqrt{1+x}-\sqrt{x})^{m}}= \frac{1}{\cosh\left(\sqrt{m^2 \, \varphi(x)} \right)}\ . \label{GLT3}
\end{equation}
\end{prop}
\vspace{5pt}
\begin{rem}
These Gauss-Laplace transforms fully characterize  the distributions of $A^{Z}_{T^{\gamma}_{c}}$ and $A^{Z}_{T^{|\gamma|}_{c}}$. By some analytic computations, we can go further into the distributional properties of these random variables. From e.g. formula \eqref{GLT2}, we can obtain the density function of $A^{Z}_{T^{|\gamma|}_{c}}$. For further details, see e.g. \cite{Vak11th,Vak11}.
\end{rem}
{\noindent \textbf{Proof of Proposition \ref{prop1}}.} Let $N$ denote a random variable following the distribution $\mathcal{N}(0,1)$.
Bougerol's identity (\ref{bougerol}) applied for $t=T^{\gamma}_{c}$ gives the identities

\begin{equation*}
\sinh(\beta_{T^{\gamma}_{c}}) \stackrel{(law)}{=} \hat{\beta}_{A^{Z}_{T^{\gamma}_{c}}}\stackrel{(law)}{=}\sqrt{A^{Z}_{T^{\gamma}_{c}}} \: N,
\end{equation*}
which in turn gives that for every fixed $c>0$,
\begin{equation}\label{invo}
\sinh(C_{c}) \stackrel{(law)}{=}\hat{\beta}_{A^{Z}_{T^{\gamma}_{c}}},
\end{equation}
where $(C_{c},c\geq0)$ is a standard Cauchy process, and we denote by $h_c$ the probability density function of $C_c$, that is
$$h_c(y)=\frac{c}{\pi(c^2+y^2)}, \quad y\in\mathds{R}.$$
Recalling by (\ref{p3}) that $\arg \sinh (y)= \sqrt{\varphi(y^2)}=\log(y+\sqrt{1+y^{2}})$ and identifying the densities of the two variables involved in (\ref{invo}),  we get
\begin{eqnarray*}
\mbox{on the LHS:} && \frac{1}{\sqrt{1+y^{2}}} \ h_{c}(\arg \sinh y) = \frac{1}{\sqrt{1+x^{2}}} \ h_{c}( \sqrt{\varphi(y^2)}) ; \\
\mbox{on the RHS:} && E \left[ \frac{1}{\sqrt{2\pi A^{Z}_{T^{\gamma}_{c}}}} \exp \left( -\frac{y^{2}}{2A^{Z}_{T^{\gamma}_{c}}} \right) \right].
\end{eqnarray*}
Performing the change the variables $x=y^2$ we get \eqref{GLT1}.  Formula (\ref{GLT2}) follows by Bougerol's identity in law applied for $T^{|\gamma|}_{c}$ and by the same arguments as previously,
since the density of $\beta_{T^{|\gamma|}_{c}}$ is given by (see e.g. \cite{BiY87})
$$\frac{1}{2c}\frac{1}{\cosh(m y)}=\frac{1}{c}\frac{1}{e^{m y}+e^{-m y}} \; .$$
\hfill \QED
\subsection{Infinite divisibility properties}\label{s2}
We shall see that formulae \eqref{GLT1} and \eqref{GLT2} yield  infinite divisibility properties for the inverse of $A^{Z}_{T^{\gamma}_{c}}$ and of $A^{Z}_{T^{|\gamma|}_{c}}$. For this purpose, we need some preparation.
\\
Let $\mathcal{BF}$ denote the class of Bernstein functions,  $\mathcal{CBF}$ the subclass of complete Bernstein functions and $\mathcal{TBF}$ the sub-subclass of Thorin Bernstein functions (see \cite{SSV} for more details). The class of infinitely divisible distributions $ID$ (resp. class of Bondesson distributions $BO$, class of Generalized Gamma Convolution distributions $GGC$) corresponds to the distribution of a positive random variable $X$ whose Laplace transform is such that
$$ E [e^{-x X}] = e^{-\phi(x)}, \;x\geq 0, \qquad  \phi \in \mathcal{BF} \; \mbox{(resp. $\mathcal{CBF}$, $\mathcal{TBF}$)}.$$
That means that $\phi$ is represented by
$$\phi(x)= d x +\int_{(0,\infty)}(1-e^{-xu})\, \nu(du)$$
where $d\geq0$, the associated L\'evy measure  $\nu$  satisfies $\int_{(0,\infty)}\min(1,u)\, \nu(du) <\infty$ and
the subclass $\mathcal{CBF}$ (resp. $\mathcal{TBF}$) corresponds to the case where $\nu$ is absolutely continuous  with density function $l$ such that
\begin{equation}\label{l}
u\mapsto l(u) \;\mbox{(resp.$\; u\mapsto ul(u)$) is completely monotone}.
\end{equation}
The inclusions $\mathcal{TBF} \subset \mathcal{CBF} \subset\mathcal{BF} \subset$  justify the fact that  $GGC \subset BO \subset ID$. Notice that  the classes $\mathcal{BF}$ and $\mathcal{CBF}$ are both convex cones stable by composition, whereas $\mathcal{TBF}$ is only a convex cone. Nevertheless, the subclass $\mathcal{TBF}$ enjoys the following property, stated as Theorem 8.4, p. 112 in \cite{SSV}: for a function $\phi \in \mathcal{TBF}$, we have that
\begin{equation}\label{sti}
\psi \circ \phi \in  \mathcal{TBF}, \; \mbox{for every} \; \psi\in \mathcal{TBF}\quad \mbox{if and only if} \quad \frac{\phi'}{\phi} \; \mbox{is a Stieltjes function},
\end{equation}
see \cite{SSV} for the definition of Stieltjes functions. Recall the function $\varphi$ introduced in (\ref{p3}). Then,
\begin{equation}
x\mapsto \varphi(x) = \arg \sinh ^2 (\sqrt{x})=\log^2(\sqrt{x}+\sqrt{1+x}) \in \mathcal{TBF}.\label{p33}
\end{equation}
\begin{rem}\label{achiev}
The latter is not observed in \cite{SSV}, but it could be obtained through the trivial equality $\sqrt{\varphi(x)}=\log(\sqrt{x}+\sqrt{x+1})=\frac{1}{2}\arg\cosh(2x+1)$
and \cite[entry 78, p. 338]{SSV}. Further, some computations stemming from \cite[entry 80, p. 338]{SSV} give that the logarithmic derivative of $\varphi$
is represented for every  $x>0$, by
\begin{equation}\label{rori}\frac{\varphi'(x)}{\varphi(x)}=\int_0^\infty e^{-xu}\, f(x) \,dx, \qquad f(x)
=e^{-\frac{x}{2}}\left(\cosh(\frac{x}{2})+\int_0^\infty I_\nu(\frac{x}{2}) d\nu\right),
\end{equation}
where $I_\nu$ is the modified Bessel function of the first kind. Showing that the logarithmic derivative of $\varphi$ is a Stieltjes function amounts
to show that the function $f$ in (\ref{rori}) is completely monotone, however this is not achieved here. Note that  if the latter is true,
then  we will have the following property that can be exploited in the representations (\ref{p1}) and (\ref{p2}) below:
\begin{equation}\label{insta}
\log(1+t\,\varphi) \in \mathcal{TBF}, \quad \mbox{\it for all }\, t>0.
\end{equation}
\end{rem}
For a positive random variable $X$, we denote by $X_{[u]}$ a version of the induced length-biased law of order $u$, that is:
\begin{equation} \label{biais}
\mbox{$X_{[u]}$ is a realization of the distribution $\displaystyle \frac{x^{u}}{E[X^u]} \, P(X\in dx)$},
\end{equation}
whenever $E[X^u]<\infty$. By \cite[Theorem 6.2.4]{Bon92}, we have that
\begin{equation}\label{bibi}
\mbox{if $X\sim GGC\quad$ and $\quad E[X^u]<\infty$ for $u<0, \quad$ then  $X_{[u]}\sim GGC$}.
\end{equation}

From now on, we adopt the following notations:  $G_{\frac{1}{2}}$ has the Gamma distribution with shape parameter $1/2$ and scale parameter 1, $\mathbf{e}_{k}$, $1\leq k\leq n$, denote $n$ independent exponentially distributed random variables with parameter 1, independent of $G_{\frac{1}{2}}$  and the length-biased random variables:
\begin{equation}\label{xci}
X_{1,c}:=\big(\frac{1}{2\, A^{Z}_{T^{\gamma}_{c}}}\big)_{[\frac{1}{2}]}\quad \mbox{and}\quad X_{2,c}:=\big(\frac{1}{2\, A^{Z}_{T^{|\gamma|}_{c}}}\big)_{[\frac{1}{2}]}.
\end{equation}
Expressions (\ref{GLT1}) and (\ref{GLT2}) are reformulated as follows:
\begin{equation}
 E_x \left[e^{- X_{1,c}} \right] = \frac{1}{\sqrt{1+x}} \, \frac{c^{2}}{c^{2}+\varphi(x)}, \qquad
 E_x \left[e^{- X_{2,c}} \right]=  \frac{1}{\sqrt{1+x}} \, f_{m}(x) , \label{T1}
\end{equation}
where $\varphi$ and $f_{m}$ are given by (\ref{p33}) and (\ref{GLT3}) respectively. The next proposition comes essentially from \cite{VaY12b} and its proof requires
the use of Chebyshev's polynomials. The second statement follows from (\ref{l}) and (\ref{bibi}).
\begin{prop}\label{propTLmint} For every integer $m$, we have the following. \\
1) The function $x\mapsto f_{m}(x)$  is the Laplace transform of a positive random variable $\mathbf{K} \sim GGC$ which has the representation:
\begin{itemize}
\item for $m=2n+1, \quad \mathbf{K}= G_{\frac{1}{2}}+\sum^{n}_{k=1} \frac{1}{a_k} \, \mathbf{e}_{k}, \qquad a_{k}= \sin^{2}\left(\frac{\pi}{2}\frac{2k-1}{2n+1}\right); \ k=1,2,\ldots,n$;
\item for $m=2n , \quad \mathbf{K}=\sum^{n}_{k=1} \frac{1}{b_k} \,\mathbf{e}_{k}, \qquad  b_{k}= \sin^{2}\left(\frac{\pi}{2}\frac{2k-1}{2n}\right); \ k=1,2,\ldots,n$.
\end{itemize}
The associated L\'{e}vy measures are:
\begin{itemize}
\item for $m=2n+1, \quad \nu (dz)= \frac{dz}{z}  \sum^{n}_{k=1} e^{-a_k\, z } $;
\item for $m=2n, \quad \nu (dz)= \frac{dz}{z}  \sum^{n}_{k=1} e^{-b_k\, z }$.
\end{itemize}
2) The random variables $X_{2,c}$ given in (\ref{xci}) satisfy the identity in law  $X_{2,c} \stackrel{(law)}{=} G'_{\frac{1}{2}} +\mathbf{K}\sim GGC$, where $G'_{\frac{1}{2}}$ is a copy of $G_{\frac{1}{2}}$, independent of $\mathbf{K}$. We also have that $1/\big(A^{Z}_{T^{\gamma}_{c}}\big)$ and $1/\big(A^{Z}_{T^{|\gamma|}_{c}}\big) \sim GGC$.
\end{prop}
{\noindent For other results and variants concerning properties of $A^{Z}_{T^{\gamma}_{c}}$ and $A^{Z}_{T^{|\gamma|}_{c}}$, the interested reader is  addressed to \cite{VaY12b,VaY12}
and to the references therein.}
\\
\bigskip

\noindent {\bf Case where $m$ is not an integer.} \smallskip

At the end of paper \cite[Remark 3.2]{VaY12b} it is conjectu that formula \eqref{GLT2} yields infinite divisibility properties for every $m>0$ (not necessarily an integer).
The next proposition proves this conjecture.
\begin{prop} For every $m>0$ and $i=1,2$, the random variable  $X_{i,c}$ given in (\ref{xci})  belongs to the class $BO$, and hence, is infinitely divisible. Moreover, as $c\to \infty$, $X_{i,c}$ converges in distribution to a Gamma distribution with scale parameter 1 and shape parameter $1/2$.
\label{propg}\end{prop}
{\noindent \textbf{Proof of Proposition \ref{propg}}.} Observe that (\ref{T1})  can be also restated, for $i=1,2$, as
\begin{equation} \label{p0}
 E \left[\exp \left( -x X_{i,c} \right) \right] = \frac{1}{\sqrt{1+x}}e^{- \psi_i(x)}=  e^{- \big(\frac{1}{2}\log(1+x)+\psi_i(x)\big)}, \qquad x\geq 0,
\end{equation}
where, with elementary computations, and using the infinite product form of the function
$$\cosh (x)=\prod_{k=1} ^\infty  (1+ \frac{x^2}{d_k}), \qquad d_k= \left(\frac{\pi}{2}(2k-1)\right)^2,$$
the functions $\psi_i,\, i=1,2$ are given by:
\begin{eqnarray}
\psi_1(x) &=& \log(1+\frac{\varphi(x)}{c^2}),\label{p1}\\
\psi_2(x)  &=& \log \cosh (\sqrt{m^2\,\varphi(x)})= \sum_{k=1} ^\infty \log\left(1+ \frac{m^2\,\varphi(x)}{d_k}\right), \label{p2}
\end{eqnarray}
where $\varphi \in \mathcal{TBF}$ is given by (\ref{p33}).  At this stage, we trivially extract the second assertion in the proposition  by letting $c\to \infty$ in (\ref{p1}) and (\ref{p2}).  It remains to show that $\psi_i$ belongs to the class $\mathcal{CBF}$.  Using the fact that $x\mapsto\log(1+x)\in \mathcal{TBF}$, expressions (\ref{p0}), (\ref{p1}), (\ref{p2}), (\ref{p33}) and the stability by composition property in $\mathcal{CBF}$,  we conclude that
$$x\mapsto \frac{1}{2}\log(1+x)+\psi_i(x)\in \mathcal{CBF}, \quad \mbox{and  hence,}\quad X_{i,c} \sim BO. $$
\hfill \QED

\smallskip

\noindent {\bf Open problem in the case where $m$ is not an integer.} \smallskip

1) Having in mind the second assertion of Proposition \ref{propTLmint} which is true for any integer $m$, we surmise that   $X_{i,c} \sim GGC$, for $i=1,2$
and for every positive number $m$. This result is not proved here because we have not obtained a closed form of the L\'evy measures of $X_{i,c}$.
To do this, one needs to check whether the functions $\psi_{i},\, i=1,2$, given by (\ref{p1}) and (\ref{p2}) belong to the class $\mathcal{TBF}$.
Then, one is tempted to check whether the function  $\varphi$ given by (\ref{p33}) satisfies (\ref{sti}) or (\ref{insta})
and this does not seem to be an easy problem to deal with, see the comments in Remark \ref{achiev}.  To summarize, if the assertion is true,
then (\ref{bibi}) will yield that
$$ \frac{1}{ A^{Z}_{T^{\gamma}_{c}}} =2\,(X_{1,c})_{[-1/2]}\sim GGC  \qquad \mbox{and}\qquad \frac{1}{ A^{Z}_{T^{|\gamma|}_c}} =2\,(X_{2,c})_{[-1/2]} \sim GGC.$$

2) If furthermore the distribution  of $X_{i,c}$ belongs to the subclass $HCM$ of $GGC$  (i.e. if the  density function of $X_{i,c}$ is  hyperbolically completely monotone,
see \cite[p. 55]{Bon92} for the  definition), then, by the comments   \cite[p. 69]{Bon92}, the random variables $A^{Z}_{T^{\gamma}_{c}}$ and $A^{Z}_{T^{|\gamma|}_{c}}$ will also have an HCM density and hence they will be infinitely divisible.

\subsection{De Blassie's result: a new proof}\label{DeBlassie}
In this section, we present a new simple proof of De Blassie's result in \cite{DeB87,DeB88}
concerning the Bessel clock $H^{Z}_{t}:=\int^{t}_{0}\frac{ds}{\left|Z_{s}\right|^{2}}$.
\begin{prop}\label{AsR}
For every $u>0$ and for every $\lambda>0$, we have that
\beqq
P\left(\int^{t}_{0}\frac{ds}{\left|Z_{s}\right|^{2}}\leq u\right)=O\left(t^{-\lambda}\right),\quad \mbox{as $t\rightarrow\infty$.}
\eeqq
\end{prop}
\noindent In order to prove Proposition \ref{AsR}, we shall make use of D. William's  "pinching method" (see e.g. \cite{Wil74,MeY82}).
Loosely speaking, when Williams studied windings of Brownian motion, instead of working out directly the asymptotics of the winding process $\theta$, he studied the asymptotic behaviour of this process taken at a random time, depending on $\theta$ (for similar results but with the use of a random time independent of $\theta$, see e.g. \cite{Vak11}). Next, one simply remarks that the difference between the initial winding process and the subordinated process is finite, and renormalising appropriately, this difference converges to 0. Hence, the asymptotic study of the renormalised subordinated process yields similar results for the renormalised initial one.
\\
{\noindent \textbf{Proof of Proposition \ref{AsR}.}}
Recall that the first passage time $T^{\beta}_{t}$, defined by (\ref{exit}), inherits the scaling property of Brownian motion as follows:
\beqq
\frac{T^{\beta}_{\sqrt{t}}}{t} \stackrel{(law)}{=}T^{\beta}_{1}.
\eeqq
Williams' "pinching method"  allows us to replace $t$ by $T^{\beta}_{\sqrt{t}}$ when $t\rightarrow\infty$. Indeed, $H^Z_{t}-H^Z_{T^{\beta}_{\sqrt{t}}}$ converges to a finite variable, as $t\rightarrow\infty$ (see also \cite{Vak11th,Vak11}).
\\
First, we choose $A,B>0$ such that
$$P\left(A<\frac{T^{\beta}_{\sqrt{t}}}{t}<B\right)=\frac{1}{2}.$$
We also remark that $T^{\beta}$ and $H^Z$ are independent, thus
\beq\label{H1}
P\left(H^Z_{t}\leq u\right)=2 P\left(A<\frac{T^{\beta}_{\sqrt{t}}}{t}<B;\;H^Z_{t}\leq u\right)=2 P\left(A\,t<T^{\beta}_{\sqrt{t}}<B\,t;\;H^Z_{t}\leq u\right).
\eeq
Moreover, since $H$ is increasing, we obtain the inequalities
\beqq
At<T^{\beta}_{\sqrt{t}}<B\,t  \;\;\Longleftrightarrow \;\;\frac{T^{\beta}_{\sqrt{t}}}{B}<t<\frac{T^{\beta}_{\sqrt{t}}}{A},
\eeqq
hence, the first part of this inequality yields
\beq\label{H2}
H^Z_{t}\leq u\Longrightarrow H^Z_{\left(T^{\beta}_{\sqrt{t}}/B\right)}\leq u.
\eeq
Now, \eqref{H1} and \eqref{H2} give
\beq\label{H3}
P\left(H^Z_{t}\leq u\right)\leq 2 P\left(H^Z_{T^{\beta}_{\sqrt{t}}/{B}}\leq u\right).
\eeq
With $a(y):= \sqrt{\varphi(y^2)}=\sinh^{-1} (y)= \arg \sinh (y)= \log\left(y+\sqrt{y^{2}+1}\right)$, $y \in \mathds{R}$, we have asymptotically that $$a(\sqrt{t}) \approx \frac{1}{2} \log t, \quad {as}\;t\rightarrow\infty.$$
If $\hat{T}^{\beta}$ is an independent copy of $T^{\beta},$ we have
\beq\label{H4}
P\left(H^Z_{t}\leq u\right) \simeq P\left(H^Z_{\hat{T}^{\beta}_{\sqrt{t}}}\leq u\right), \quad \mbox{as} \;t\rightarrow\infty.
\eeq
Following \cite{Vak11th,Vak11}, the skew-product representation of planar Brownian motion \eqref{skew-product} and Bougerol's identity in law \eqref{bougerol} yield that
if $(\delta_t, t\geq0)$ denote another independent real Brownian motion, then, with obvious notation, we have that (see also \cite[Proposition 2.3]{Vak11}).
\beq\label{H5}
H^Z_{T^{\delta}_{b}} \stackrel{(law)}{=} T^{\beta}_{a(b)}, \quad \mbox{for every $b\geq0$}.
\eeq
Indeed, using the symmetry principle (\cite{And87} for the original note and \cite{Gal08} for a detailed discussion), Bougerol's
identity in law \eqref{bougerol} is equivalently stated as:
\beqq
    \sinh(\bar{\beta}_{u}) \stackrel{(law)}{=} \bar{\delta}_{A^{Z}_{u}(\beta)}, \quad \mbox{for any fixed} u>0.
\eeqq
Hence, identifying the densities of the two parts and recalling that $H^{Z}$ is given by (\ref{clock}), we get easily \eqref{H5}.
\\ \\
\noindent Now, using \eqref{H5} with $b=\sqrt{t}$ and adapting appropriately the notation, \eqref{H4} writes
\beqq
P\left(H^Z_{t}\leq u\right)=P\left(T^{\beta}_{a(\sqrt{t})}\leq u\right) P\left(T^{\beta}_{\frac{1}{2} \log t}\leq u\right), \quad \mbox{as} \; t\rightarrow\infty.
\eeqq
Below, the symbol "$\approx$" means "is of order". We have
\beqq
P\left(T^{\beta}_{h}\leq u\right)\approx\frac{\sqrt{u}}{h} \exp\left(-\frac{h^{2}}{2u}\right)\stackrel{h=\frac{1}{2} \log t}{\approx}
\frac{2\sqrt{u}}{\log t} \exp\left(-\frac{(\log t)^{2}}{8u}\right).
\eeqq
Thus, with $h=\frac{1}{2}\log t$, which corresponds to the asymptotic behaviour for $t\rightarrow\infty$, we get that
\beq\label{sigma}
P\left(T^{\beta}_{h}\leq u\right)\leq\frac{2\sqrt{u}}{\log t} \exp\left(-\frac{(\log t)^{2}}{8u}\right).
\eeq
Observe that for big values of $d$, we have $\exp(-d^{2})\leq\exp(-\lambda d), \ \forall \lambda>0$, hence, with $d=\log t$, we get
$$\exp\left(-(\log t)^{2}\right)\leq \exp(-\lambda(\log t))=\frac{1}{t^\lambda} , \qquad \forall \lambda>0.$$
Using \eqref{H3}, \eqref{sigma} and last elementary remark, we obtain
\beqq
P\left(H^Z_{t}\leq u\right)\leq  \frac{C_{u,B,\lambda}}{t^{\lambda}}, \qquad \forall \lambda>0,
\eeqq
where $C_{u,B,\lambda}$ denotes a positive constant depending on $u,B$ and $\lambda$. This finishes the proof.
\hfill \QED

\section{Planar Stable Processes}\label{St}

\subsection{The winding process}\label{subStable}
In this section, we  focus on isotropic planar Stable processes.
Bertoin and Werner \cite{GVA86} obtained the following results for $\alpha\in(0,2)$ (see \cite{BeW96} for the proofs ). Let us denote now by $dz$ the Lebesgue measure on $\mathbb{C}$ and for every complex number $z\neq0$, $\omega(z)$ stands for the determination of its argument valued in $\left(\right.-\pi,\pi\left.\right]$.
\begin{lem}\label{lemma1}
The time-changed process $(\theta^{U}_{A^{U}(u)},u\geq0)$ is a real-valued symmetric L\'{e}vy process, say $\rho$. It has no Gaussian component and its L\'{e}vy measure has support in $[-\pi,\pi]$.
Moreover, the L\'{e}vy measure of $\theta_{A^{U}(\cdot)}$ is the image of the L\'{e}vy measure of $U$ by the mapping $z\rightarrow\omega(1+z)$. Consequently, $E[(\theta_{A^{U}(u)})^{2}]=u\, k(\alpha)$, where
\beqq
k(\alpha) = \frac{\alpha \ 2^{-1+\alpha/2} \Gamma(1+\alpha/2)}{\pi\Gamma(1-\alpha/2)} \int_{\mathbb{C}} |z|^{-2-\alpha} |\omega(1+z)|^{2} dz \ .
\eeqq
\end{lem}
{\noindent For the process $U$, we use the analogue of the skew product representation for planar BM
which is the Lamperti correspondence for stable processes.
Hence, there exist two real-valued L\'{e}vy processes $(\xi_{u},u\geq0)$ and $(\rho_{u},u\geq0)$, the first one is non-symmetric
whereas the second one is symmetric, both starting from 0, such that
$$\log\left|U_{t}\right|+i\theta^{U}_{t}=\left(\xi_{u}+i\rho_{u}\right)\Bigm|_{u=H^{U}_{t}} \ .$$ }
\begin{rem} The processes $|Z|$ and $Z_{A^{U}(\cdot)}/|Z_{A^{U}(\cdot)}|$ are not independent. This is easily seen since
$|Z_{A^{U}(\cdot)}|$ and $Z_{A^{U}(\cdot)}/|Z_{A^{U}(\cdot)}|$ jump at the same times, hence they cannot be independent.
Moreover, $A^{U}(\cdot)$ depends only upon $|Z|$, hence $|Z|$ and $Z_{A^{U}(\cdot)}/|Z_{A^{U}(\cdot)}|$ are not independent.
For further discussion on the independence, see e.g. \cite{LiW11}, where it is shown that an isotropic
$\alpha$-self-similar Markov process has a skew-product structure if and only if its radial and its
angular part do not jump at the same time.
\end{rem}

\subsection{The asymptotic behaviour of windings}
{\noindent Bertoin and Werner in \cite{BeW96} obtained an asymptotic result, which is in some sense, a version of Spitzer's asymptotic Theorem \ref{Spitzertheo} for isotropic Stable L\'{e}vy processes of index $\alpha\in(0,2)$. }
\begin{theo}\label{SpiBW}
As $c\rightarrow\infty$, the family of processes  $\,\left(c^{-1/2}\theta^{U}_{\exp(ct)},t\geq0\right)\,$ converges in distribution on the space $D(\left[\right.0,\infty\left.\right),\mathbb{R})$ endowed with the Skorohod topology, to $\left(\sqrt{r(\alpha)}\beta_{t},t\geq0\right)$, where $\left(\beta_{s},s\geq0\right)$ is a real valued Brownian motion and
\beq\label{cst}
r(\alpha)=\frac{\alpha \ 2^{-1-\alpha/2}}{\pi} \int_{\mathbb{C}} |z|^{-2-\alpha} |\omega(1+z)|^{2} dz \ .
\eeq
\end{theo}
\noindent{\textbf{Proof of Theorem \ref{SpiBW}.}}
We refer to two different proofs:
\begin{enumerate}[(i)]
\item Bertoin and Werner (1996) \cite{BeW96}, using an "Ornstein-Uhlenbeck type" process and ergodicity arguments, and
\item Doney and Vakeroudis (2012) \cite{DoV12}, using the continuity of the composition function $\rho_{H^{U}(\cdot)}$ (see \cite{Whi80}).
\end{enumerate}
\hfill\QED

\section{Applications to the pricing of Asian options}\label{Asian}

\subsection{Asian options and exponential functionals of Brownian motion}
In this subsection, we return to the initial financial mathematics problem, that is, the characterisation of the distribution of
$$A^{Z}_t=\int_{0}^{t}\exp(2\beta_u)du,$$
in order to compute $E\left[\left(\frac{1}{t}A^{Z}_t-K\right)^{+}\right]$.
To that end, one may use the previously stated results to access the distribution of $A_t$ via William's o called "pinching method" \cite{Wil74,MeY82} that was also used in Subsection \ref{DeBlassie}.
We propose here to mimic again this method for our benefit, by invoking the time changes discussed in the previous sections.
\begin{prop}\label{propnew}
The following convergence in law holds
\begin{equation*}
\frac{1}{t} \log A^{Z}_{t^2}\overset{{(law)}}{\underset{t\rightarrow\infty}\longrightarrow} 2|\beta|_{T^{\gamma}_{1}}\stackrel{(law)}{=}2|C|_1,
\end{equation*}
where $C_1$ is a standard Cauchy random variable.
\end{prop}
{\noindent \textbf{Proof:}}
First, observe that
$$\log \left(\frac{A^{Z}_{T^{\gamma}_{t}}}{A^{Z}_{t^2}}\right)=\log \left(\frac{\int_{0}^{T^{\gamma}_{t}}\exp(2\beta_u)du}{\int_{0}^{t^2}\exp(2\beta_u)du}\right),$$
which is a random variable that exists (and which seems to be of no other interest here).
Renormalising by $t$, we get
\begin{equation*}
   \frac{1}{t}  \left( \log A^{Z}_{T^{\gamma}_{t}} - \log A^{Z}_{t^2}\right)= \frac{1}{t} \log \left(\frac{A^{Z}_{T^{\gamma}_{t}}}{A^{Z}_{t^2}}\right) \overset{{(law)}}{\underset{t\rightarrow\infty}\longrightarrow} 0.
\end{equation*}
Hence, studying asymptotically $t^{-1}\log A^{Z}_{T^{\gamma}_{t}}$, as $t\rightarrow\infty$, would yield similar results for $t^{-1}\log A^{Z}_{t^2}$.
Following \cite{VaY12}, applying the scaling property of Brownian motion and making a change of variables, we have that
$$A^{Z}_{T^{\gamma}_{t}}=\int_{0}^{T^{\gamma}_{t}}e^{2\beta_v}dv\stackrel{(law)}{=}t^2\int_{0}^{T^{\gamma}_{1}}e^{2t\beta_u}du, \quad \mbox{(recall that $T^{\gamma}_{t}\stackrel{(law)}{=}t^2 T^{\gamma}_{1}$ )},$$
so that, for all $t>0$, we have that
$$\frac{1}{t} \log A^{Z}_{T^{\gamma}_{t}} \,\stackrel{(law)}{=} \,\frac{1}{t} \log \left(t^2\int_{0}^{T^{\gamma}_{1}}e^{2t\beta_u}du\right) =\frac{2 \log t}{t}+\log \left(\int_{0}^{T^{\gamma}_{1}}e^{2t\beta_u}du\right)^{1/t}.$$
Using the fact that the $p$-norm converges to the $\infty$-norm when $p\rightarrow\infty$,
the latter converges for $t\rightarrow\infty$ towards $\,2\sup_{0\leq u\leq T^{\gamma}_{1}} \beta_u.\,$
By the reflexion principle (see e.g. \cite{ReY99}), we have
$$\sup_{0\leq u\leq T^{\gamma}_{1}} \beta_u\stackrel{(law)}{=}|\beta|_{T^{\gamma}_{1}}\stackrel{(law)}{=}|C_1|$$
and we deduce that
$$\frac{1}{t} \log A^{Z}_{T^{\gamma}_{t}}\overset{{(law)}}{\underset{t\rightarrow\infty}\longrightarrow} 2|C_1| \ .$$
The result for $A^{Z}_{t}$ follows immediately.
\hfill \QED
\\ \\
The distribution of $A^{Z}_{t}$ may also be characterized by a result due to Dufresne \cite{Duf00} that we state now.
For sake of completeness we shall also sketch the proof.
\begin{prop}
For every $t>0$, $x\geq 0$ and with $\varphi$ given in (\ref{p33}), we have that
\beqq
E \left[ \frac{1}{\sqrt{2\pi A^{Z}_{t}}} \exp \left( -\frac{x}{2A^{Z}_{t}} \right) \right] = \frac{1}{\sqrt{2\pi t}} \: \frac{1}{\sqrt{1+x}} \: \exp \left( -\frac{\varphi(x)}{2t} \right).
\eeqq
\end{prop}
{\noindent \textbf{Proof:}}
We appeal again to Bougerol's identity in law: for every $t>0$ fixed,
\beqq
\sinh(\beta_{t}) \stackrel{(law)}{=} \hat{\beta}_{A^{Z}_{t}(\beta)},
\eeqq
and we identify the densities of the two parts, i.e.
\begin{eqnarray*}
\mbox{on the LHS:} && \frac{1}{\sqrt{2\pi t}} \: \frac{1}{\sqrt{1+y^{2}}} \: \exp \left( -\frac{\varphi(y^2)}{2t} \right) , \\
\mbox{on the RHS:} && E \left[ \frac{1}{\sqrt{2\pi A^{Z}_{t}}} \exp \left( -\frac{y^{2}}{2A^{Z}_{t}} \right) \right].
\end{eqnarray*}
The proof finishes by the change of  variables $x=y^{2}$.
\hfill \QED

\begin{cor}For every $t>0$, we have that $1/A^{Z}_{t} \sim GGC$.
\end{cor}
{\noindent \textbf{Proof:}} With the notations (\ref{biais}),  observe that
$$E \left[\exp \left( -\frac{x}{2} \left(\frac{1}{A^Z_t}\right)_{[\frac{1}{2}]}  \right) \right] =   e^{- \chi(x)},\quad \mbox{where}\quad x \mapsto \chi(x)=  \frac{1}{2}\log(1+x)+\frac{\varphi(x)}{2t}   \in \mathcal{TBF},$$
and conclude as in the proof of Proposition \ref{propg}.
\hfill \QED
\begin{rem}
These results may easily be generalized for the functional
$$A^{(\nu)}_{t}=\int^{t}_{0}  \exp(\beta_{s}+\nu s) ds.$$
Indeed, we have access to its distribution by the following relation (see e.g. \cite{ADY97} or \cite{Vak12}): \\
with $\nu,\mu$ two real numbers, for every $t>0$ fixed ($\beta$, $B$ and $\delta$ are three independent Brownian motions),
\begin{equation*}
\sinh (Y^{(\nu,\mu)}_{t})\stackrel{(law)}{=} \int^{t}_{0} \exp(\beta_{s}+\nu s)d(B_{s}+\mu s)=\delta_{\int^{t}_{0} \exp\left(2(\beta_{s}+\nu s)\right)ds},
\end{equation*}
where $(Y^{(\nu,\mu)}_{t},t\geq0)$ is a diffusion with infinitesimal generator
\beqq
\frac{1}{2} \ \frac{d^{2}}{dy^{2}}+ \left(\nu \tanh(y) + \frac{\mu}{\cosh(y)}\right) \frac{d}{dy} \ ,
\eeqq
starting from $y= \arg \sinh (x)$.
Here, without loss of generality we may consider $\mu=0$ and mimic the approach where also $\nu=0$ which was presented above.
\end{rem}

\subsection{Asian options and exponential functionals of L\'{e}vy processes}
In this subsection, we discuss the case of Asian options in relation with L\'{e}vy processes, that is the case
where the exponential functional of interest is
$$A^{U}_{t}:= \int^{t}_{0}\exp(\alpha\xi_{s})ds.$$
Recall from Subsection \ref{subStable} that $U$ is an isotropic planar Stable process, and $\xi$, $\rho$
are two real-valued L\'{e}vy processes, the first one is non-symmetric and the second one is symmetric.
Following \cite[Subsection 6.3]{BeΥ05} and \cite{Yor01}, this case may be conside as a natural generalization of the case of Asian options where the exponential functional is associated to a Brownian motion. More precisely, as one can see in the above references, the computation of the price of Asian options corresponds to the study of the law of the exponential functional associated to a L\'{e}vy process $\xi$ at a fixed time $t$. In particular, the problem may be uced by substituting $t$ by an exponential (random) time. Hence, we have:
\begin{equation}\label{relATstable}
T^{\theta^{U}}_{c}=\inf\{ t:\theta^{U}_{t}=c \}
=(H^{U})^{-1}_{u}\Bigm|_{u=T^{\rho}_{c}}=\int^{T^{\rho}_{c}}_{0}ds\exp(\alpha\xi_{s})=: A^{U}_{T^{\rho}_{c}} \ ,
\end{equation}
and similarly
\begin{equation*}
T^{|\theta^{U}|}_{c}=A^{U}_{T^{|\rho|}_{c}} \ .
\end{equation*}
We state the following proposition only for $A^{U}_{T^{\rho}_{c}}$, a similar result may also  be obtained for $A^{U}_{T^{|\rho|}_{c}}$.
\begin{prop}
The following convergence in law holds
\begin{equation}\label{relATstable1}
\frac{1}{t} \log A^{U}_{t^{\alpha/2}}\overset{{(law)}}{\underset{t\rightarrow\infty}\longrightarrow} T^{\beta}_{\sqrt{1/r(\alpha)}} \ ,
\end{equation}
where $r(\alpha)$ is given by \eqref{cst}, $\beta$ denotes  a real Brownian motion and ${(T^{\beta}_{u})}_{u>0}$,    given by (\ref{exit})  is a $\frac{1}{2}-$stable subordinator.
\end{prop}
{\noindent \textbf{Proof:}}
Mimicking the approach of the previous subsection we can extend this result to $A^{U}_{t}$. Indeed, we easily show that
$A^{U}_{T^{\rho}_{\sqrt{t}}}-A^{U}_{t^{\alpha/2}}$ is a variable that exists, hence
$$\frac{1}{t} \left(\log A^{U}_{T^{\rho}_{\sqrt{t}}}-\log A^{U}_{t^{\alpha/2}}\right)=\frac{1}{t} \left(\log \frac{A^{U}_{T^{\rho}_{\sqrt{t}}}}{ A^{U}_{t^{\alpha/2}}}\right)
\overset{{(law)}}{\underset{t\rightarrow\infty}\longrightarrow} 0.$$
Now, following \cite{DoV12}, we use \eqref{relATstable} and Theorem \ref{SpiBW} in order to get
\begin{eqnarray*}
\frac{1}{t} \log A^{U}_{T^{\rho}_{\sqrt{t}}}&=& \frac{1}{t}\log\left(T^{\theta^{U}}_{\sqrt{t}}\right)= \frac{1}{t} \log \left(\inf\left\{u: \theta^{U}_{u}>\sqrt{t}\right\}\right) \nonumber \\
&\stackrel{u=\exp(ts)}{=}& \frac{1}{t} \log \left(\inf\left\{e^{ts}:\frac{1}{\sqrt{t}} \ \theta^{U}_{\exp(ts)}>1\right\}\right) \nonumber \\
&=& \inf\left\{s:\frac{1}{\sqrt{t}} \ \theta^{U}_{\exp(ts)}>1\right\} \nonumber \\
&\overset{{(law)}}{\underset{t\rightarrow\infty}\longrightarrow}& \inf\left\{s:\beta_{r(\alpha)s}>1\right\} \nonumber \\
&=& \inf\left\{s:\sqrt{r(\alpha)}\beta_{s}>1\right\}=: T^{\beta}_{\sqrt{1/r(\alpha)}} \ ,
\end{eqnarray*}
which finishes the proof.
\hfill \QED

\begin{cor}\label{cor}
Let $N\sim\mathcal{N}(0,1)$, The following convergence in law holds:
  \begin{equation*}
(A^{U}_{t^{\alpha/2}})^{1/t} \;\overset{{(law)}}{\underset{t\rightarrow\infty}\longrightarrow} \;e^{\frac{1}{r(\alpha)}\frac{1}{N^2}}\sim GGC.
\end{equation*}
\end{cor}
{\noindent \textbf{Proof:}}
Just use \eqref{relATstable1}, the scaling property of Stable processes, the fact that $1/N^2 \stackrel{(law)}{=} T^\beta_1   \sim GGC$ and \cite[Theorem 3]{Bon15}.
\hfill \QED
\begin{rem}
The result of Corollary \ref{cor} could be useful in order to obtain asymptotic closed formulae about "jump type" Asian option prices by following the spirit of Geman and Yor \cite{GeY01}. This problem will be further discussed in a forthcoming paper.
\end{rem}

\vspace{20pt}
\noindent\textbf{Acknowledgements} \\
The  first author  would like to extend his sincere appreciation to the Deanship of Scientific Research at King Saud University for funding this Research group No. (RG-1437-020). The research of the second author was partly financed by the Project: Postdoctoral Researchers 2016-2017 of the University of Cyprus. He is indebted to Professor Konstantinos Fokianos (University of Cyprus - UCY, Department of Mathematics and Statistics) for supervising his postdoctoral stay at the University of Cyprus where he prepa several parts of this work. Moreover, the authors would like to thank anonymous referees for useful comments that improved the paper and in particular Subsection \ref{s2}. Finally, they would like to express once more their gratitude to Professor Marc Yor who passed away suddenly some years ago. The stimulating discussions they had with him remain always a source of inspiration.



\end{document}